\documentclass[12pt]{amsart}
\usepackage{amscd}
\usepackage{verbatim}
\usepackage{amssymb, amsmath, amsthm, amscd,ifthen}
\usepackage[dvips]{graphics}
\usepackage[cp866]{inputenc}
\usepackage{graphicx}
\usepackage{epsfig}

\usepackage{amsmath,amssymb,amscd,}
\usepackage[all]{xy}

\usepackage{amssymb}

\newtheorem{thm}{Theorem}[section]

\newtheorem{Ex}[thm]{Example}

\newtheorem{lemma}[thm]{Lemma}

\theoremstyle{definition}

\title[Morse index of a cyclic polygon]
{ Morse index of a cyclic polygon}

\author[G. Panina, A. Zhukova]
       {Gaiane Panina, Alena Zhukova}


\begin{document}
\begin{abstract}
It is known that cyclic configurations of a planar polygonal linkage
are  critical points of the signed area function.  In the paper, we
announce an explicit formula of the Morse index for the signed area
of a cyclic configuration.

It depends not only on the combinatorics of a cyclic configuration,
but also includes some metric characterization.
\end{abstract}
\maketitle

\setcounter{section}{0}
\section{Introduction}
We study planar polygonal linkages with two pinned vertices and
their flexes in the plane with allowed self-intersections.

Let us make this precise.

A sequence of real positive numbers $L=(l_1,...,l_n)$ which can be
realized as a sequence of edge lengths of a planar polygonal chain
is called an $n$-\textit{linkage}. A linkage carries a natural
orientation which we indicate in figures by arrows.

The sequence of points in $P=(p_1,...,p_n), \ p_i \in \mathbb{R}^2$
is called a \textit{configuration of the linkage}  $L$, if
\begin{enumerate}
    \item $l_i=|p_i,p_{i+1}|$, i.e. the lengths are fixed, and
    \item $p_1=(0,0); \ p_2=(0,l_1)$, i.e., the first two vertices
    are pinned  down.
\end{enumerate}

A configuration $P$  is called \textit{cyclic} if all its vertices
lie on a circle.

By $\mathcal{M}(L)$  we denote the \textit{moduli space of} $L$,
i.e., the set of all configurations of $L$. Generically,
$\mathcal{M}(L)$ is a smooth manifold of dimension $n-3$. It embeds
canonically in $ \mathbb{R}^{2(n-2)}$  by listing the coordinates of
all the vertices of a configuration except for the first two ones
(which are pinned down).

The core object of the paper is  the \textit{signed area} $A(P)$ of
a configuration $P$. For a generic linkage $L$, the signed area
$A(P)$ is a Morse function on $\mathcal{M}(L)$.

It was known since long that $A$ achieves its maximum at the convex
positively oriented cyclic configuration of $L$. Consequently, $A$
achieves its minimum at the \textit{anticonvex} (convex negatively
oriented) cyclic configuration of $L$.

 The interpretation of cyclic configurations as
critical points of the signed area function was suggested in
\cite{khi2}. As was shown in \cite{eljikh}, this is indeed the case
for generic planar quadrilaterals and pentagons. The result was
extended in \cite{khipan} by proving that the same holds for generic
cyclic configurations   with arbitrary number $n$ of vertices.

The aim of the paper is to compute the Morse index $m(P)$ for a
cyclic configuration $P$.  This will be done  in Theorems
 \ref{signHess} and \ref{mainthm_sequence}.

To our opinion, the formula obtained leaves much room for other
interpretation, equivalent reformulations and further study.

The paper only announces the results. That is, we omit details of
the proofs and the delicate discussion about \textit{genericity }of
 linkages, cyclic configurations and deformations.

\section{Preliminaries}

\begin{thm} \cite{khipan} Let $L$ be a
a  generic linkage. Its  configuration $P$ is a critical point of
the function $A$  iff $P$ is a cyclic configuration.\qed
\end{thm}

\bigskip
 Given a
cyclic configuration $P$ of a linkage $L$, we use throughout the
paper the following notation:

$r_P$ is the radius of the circle  which circumscribes $P$.

$O$ is the center of the circle.

$\{p_i\}_{i=1}^n$ are the vertices of $P$. We assume that the
numeration is cyclic, that is, for instance, $p_0=p_n, \
p_{n+1}=p_1$.

$l_i$  is the length of the $i$-th edge.

$\alpha_i$  is the half of the angle between the vectors
$\overrightarrow{Op_i}$ and $\overrightarrow{Op_{i+1}}$. The angle
is defined to be positive, orientation is not involved.

$m(P)$ is the Morse index of the signed area function $A$.

 $Hess(P)=D^2A$  is the
Hessian of $A$ at the point $P$.

$H(P)=Det(Hess(P))$  is the determinant of the Hessian.

$\mathcal{H}(P)$ is the sign of $H(P)$.

\begin{figure}\label{firstExample}
\centering
\includegraphics[width=8 cm]{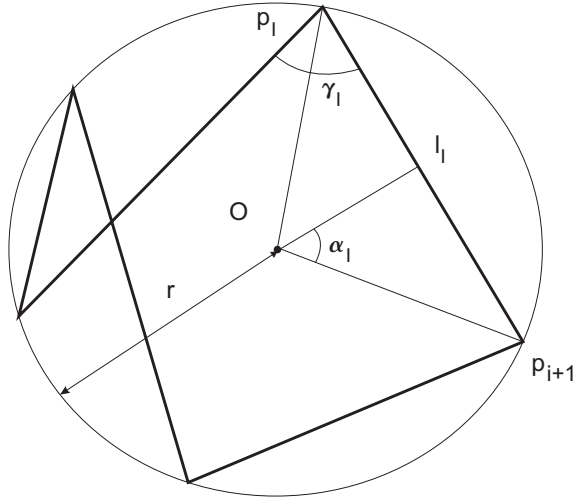}
\caption{Basic notation }
\end{figure}

A cyclic configuration $P$ is called \textit{central} if one of its
edges passes through the center of the circle.

For a non-central configuration, define

$\varepsilon_i$ is the orientation of the edge $p_ip_{i+1}$, that
is,

 $\varepsilon_i=\left\{
                       \begin{array}{ll}
                         1, & \hbox{if the center $O$ lies to the left of } p_ip_{i+1};\\
                         -1, & \hbox{if the center $O$ lies to the right of } p_ip_{i+1}.
                       \end{array}
                     \right.$

$E(P)$ is the string of orientations of all the edges, that is,
$E(P)=(\varepsilon_1,...,\varepsilon_n)$

For a non-central configuration,  put

$$\delta P=\Sigma \varepsilon_i \tan \alpha_i, \ \ \hbox{and}$$

$d(P)$ is the sign of $\delta P$.

\begin{figure}\label{firstExample}
\centering
\includegraphics[width=6 cm]{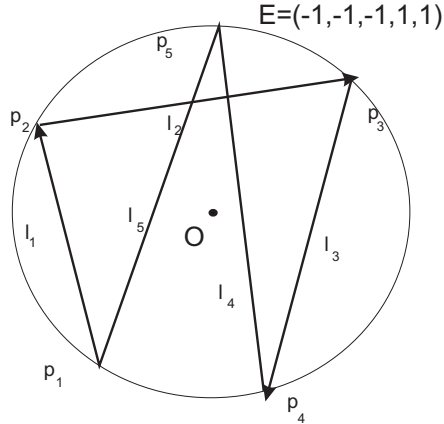}
\caption{The value of $E$ }
\end{figure}

Each cyclic configuration $P$  is uniquely defined by the pair
$(r_P;E(P))$.

\bigskip

We start with two examples. The first one gives us the base for
computation of $m(P)$. The other one is just an illustration.

\begin{Ex}\cite{eljikh} For a  generic 4-gonal  linkage, there are two
possible cases:
\begin{enumerate}
    \item The configuration   space $\mathcal{M}(L)$ is disconnected.
    Then $L$ has four cyclic configurations listed in Fig. 3. The Morse index depends on
the sign of the area $A$ and on the self-intersection of the
configuration.
    \item The  configuration  space $\mathcal{M}(L)$ is connected.
    Then $L$ has  two cyclic configurations (the first two ones  in Fig. 3).
\end{enumerate}
\end{Ex}

\begin{figure}\label{4-gon}
\centering
\includegraphics[width=10 cm]{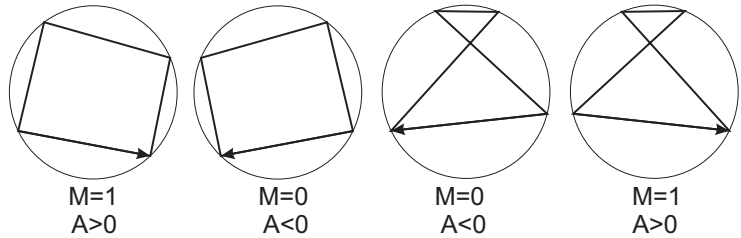}
\caption{4-gonal cyclic configurations}
\end{figure}

\newpage

\begin{Ex} For  the equilateral pentagonal linkage $L=(1,1,1,1,1)$ there are 14 cyclic
configurations  listed in the Fig.  4.

(1). The convex and the anticonvex ones are the global maxima of the
signed area $A$ (their Morse indices are $2$ and $0$ respectively).

(2). The starlike   configurations are local maximum and a local
minimum of $A$.

(3).  There are $10$ more configurations that have three consecutive
edges aligned. Their  Morse indices equal $1$.

This can be either deduced from Theorem \ref{mainthm_sequence}, or
obtained by simple symmetry reasons.
\end{Ex}

\begin{figure}\label{pentagon}
\centering
\includegraphics[width=10 cm]{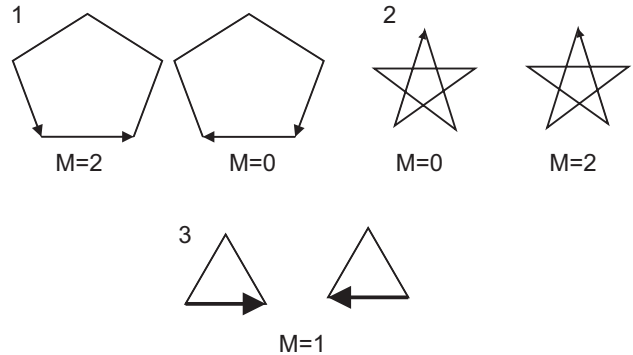}
\caption{Cyclic configurations of the equilateral pentagonal linkage
and their Morse indices }
\end{figure}

\section{Dynamics of Morse points. Computation of $\mathcal{H}$. }

A \textit{deformation of a linkage}  is a  one-parametric continuous
family $L(t),  \  t \in [0,1] $  of linkages.

In particular, we will explore
 \textit{cyclic deformations of a linkage} $L$  which arise
through the following construction.

Given a cyclic configuration $P$ of a linkage $L$, we fix the radius
$r_P$ and force the vertices $p_i$ move along the circle. We get a
 continuous family $L(t), \ L(0)=L$  together with a
 continuous family of its cyclic configurations $P(t), \ P(0)=P$.

It can happen that during such a deformation, two consecutive
vertices $p_i$ and $p_{i+1}$   meet, the edge $l_i$ vanishes and
then appears again with a different orientation $\varepsilon_i$ (see
Figure 5). This will be called \textit{a flip}.

Given a cyclic configuration $P$, we will apply  a generic cyclic
deformation $P(t)$ and treat all the moments $t$ when $\mathcal{H}$
changes.

Here is a rough idea to be explored below:  during a  deformation
$L(t)$, the Morse points move and sometimes several of them meet. By
Cerf theory \cite{cerf}, if a Morse point meets no other Morse
point, its Morse index and the value of  $\mathcal{H}$   do not
change.  If two Morse points meet,  their Morse indices satisfy
$|m_1-m_2|=1$. Consequently, the values of $\mathcal{H}$ at these
points are different.

Here is the precise construction:

Let $L=L(t), \ \ t \in [0,1]$ be a cyclic deformation of a generic
linkage $L$.

\begin{figure}\label{flip}
\centering
\includegraphics[width=10 cm]{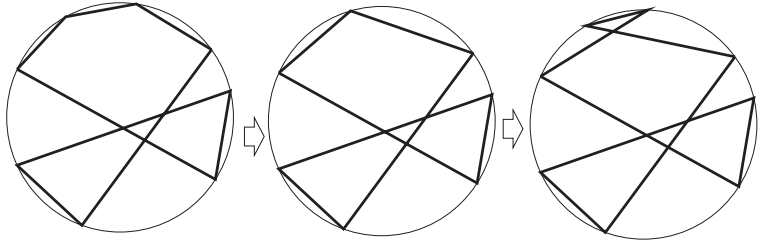}
\caption{A cyclic deformation passes through a flip }
\end{figure}

By generity reasons we can assume that
     all of $P(t)$, except for a finite number of $t$, are Morse
    points. Besides, we assume that no two of the following events
    happen at the same moment $t$.

\begin{enumerate}
    \item $P(t)$ is a central configuration (one of the edges passes
    through $0$).
    \item One of the edges vanishes ($t$ is the moment of a flip).
    \item $\delta(P(t))=0$.
\end{enumerate}

\begin{lemma}\label{thm-delta-h}
Assume that for some $t_0 \in [0,1]$, the configuration $P(t_0)$ is
neither central nor has a vanishing edge. Then generically,
\begin{enumerate}
    \item $\mathcal{H}(P(t_0))=0$ iff  $\mathcal{H}(P(t))$ changes at
    $t=t_0$.
    \item $\delta P(t_0)=0$ iff  $d(P(t))$ changes at
    $t=t_0$.
    \item $\mathcal{H}(P(t))$ changes at $t=t_0$ if and only if
$\delta P(t)$ changes its sign  at $t=t_0$.
\end{enumerate}

\end{lemma}
Sketch of the proof.

(1) and (2) follow from analysis of the formulae for $\mathcal{H}$
and $\delta P$.

 Before proving (3), let us first make the following observation. Assume
that a linkage $L$ and a string
$E=(\varepsilon_1,...,\varepsilon_n)$ are fixed. There exists a
cyclic configuration $P$  of $L$ with $E(P)=E$ inscribed in a circle
of radius $r_P$ if and only if for some integer $k$, the number
$r_P$ is a root of the function  $F_{L,E}$, where
$$F_{L, E}(r)=\sum_{i=1}^n \varepsilon_i \arcsin(
\frac{l_i}{2r})-\pi k  .$$

Now prove (3).

 1.
 Assume that  $\delta P(t_0)=0$. Then  the
function  $F_{L, E}$  has a multiple root at $r_P$. Indeed, it is a
root just by the above observation. Besides, we easily have
$$-\frac{dF_{L(P(t_0)),\ E}}{dr}|_{r=r_P} =\sum_{i=1}^n \frac{\varepsilon_i
l_i}{2r_P^2\sqrt{1-\sin^2 \alpha_i}} =
\frac{1}{2r_P}\sum_{i=1}^n\varepsilon_i \tan
\alpha_i=\frac{1}{2r_P}\delta P(t_0)=0$$

(Here we write for short $l_i$ and $\alpha_i $  instead of
$l_i(P(t_0))$ and $\alpha_i (P(t_0))$.)

This means that the linkage  $L(P(t_0+\varepsilon))$  has two Morse
points $P(t_0+\varepsilon)$ and $P'(t_0+\varepsilon)$ with  one and
the same $E(P)=E(P')=E$
   such that  $P(t_0+\varepsilon)$ and $P'(t_0+\varepsilon)$ tend to $P(t_0)$ as $\varepsilon$ tends to zero.

 By  above
 arguments of Cerf theory and continuity reasons,
$\mathcal{H}(P(t_0))=0$.

2. Conversely, if  $\mathcal{H}(t_0)=0$, then $P$ is a limit point
of two critical points of a deformation $L(q)$. Hence $r_P$ is a
multiple root of $F_{L, E}$.\qed

\begin{lemma}\label{lemma-flip}Let a generic cyclic deformation $P(t)$ pass  through a flip. Then the values of
$\mathcal{H}(P)$ and $E(P)$ change, whereas the value $d(P)$ does
not change.\qed
\end{lemma}

\begin{lemma}  If during a generic cyclic deformation $P(t)$ passes a central
configuration, then  $d (P)$ and $E(P)$ change, but $\mathcal{H}(P)$
remains constant.\qed
\end{lemma}

Taken together, the three lemmata  imply the first main result of
the paper:

\begin{thm}\label{signHess} Let $P$ be a generic cyclic configuration.

Denote  by $e(P)$ the number of positive entries in
$E(P)=(\varepsilon_1,...,\varepsilon_n)$.

Remind that  $d(P)$ the sign of  $\delta(P)$. Then

$$\mathcal{H}(P)= -d(P) (-1)^{e(P)} $$

\end{thm}

Proof. Consider a generic cyclic deformation  joining $P$ and some
convex cyclic configuration $P_{CONV}$ with the same number of
vertices.

For $P_{CONV}$ the theorem is obviously valid. Besides, the above
three lemmata imply that the product  $\mathcal{H}(P(t))\cdot
d(P(t))(-1)^{e(P(t))} $ does not change as $t$ changes. \qed

\section{Computation of the  Morse index }
The following lemma provides an inductive computation of  $m(P)$.
\begin{lemma}\label{decomp}Let $P$ be a generic cyclic configuration. Adding its diagonal $p_ip_j$,
we express $P$ as the homological  sum of a cyclic configurations
$P_1$ and a trigonal cyclic configuration $P_2$  (see Fig. 6). Then
\begin{enumerate}
    \item Either  $m(P)=m(P_1)$   or $m(P)=m(P_1)+1.$
    \item The Morse index   $m(P_1)$ together with $\mathcal{H}(P)$  uniquely
    determine $m(P)$.
\end{enumerate}

\end{lemma}
Sketch of the proof:

(1) follows from the fact that  the moduli space $\mathcal{M}(P_1)$
is a codimension one submanifold  of $\mathcal{M}(P)$; (2) easily
follows from (1).\qed

\begin{figure}\label{homological-sum}
\centering
\includegraphics[width=8 cm]{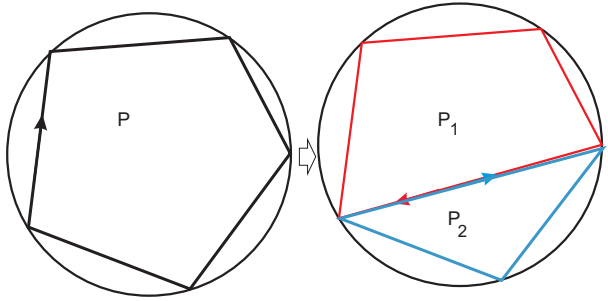}
\caption{Homological decomposition $P=P_1+ P_2$ }
\end{figure}

\newpage

Iterative application of the above lemma starting from 4-gonal
linkages immediately implies the final theorem:

\begin{thm}\label{mainthm_sequence}
Let $P=(p_1,...,p_n)$ be a generic cyclic configuration. Introduce
its subconfigurations  $P_3,...,P_n$  (see Fig. 7)
$$P_i=(p_1,...,p_i), \ \ i=3,...,n.$$

 Put
 $\mathcal{H}(P_3)=1$ for the trigonal configuration $P_3$.

Then the Morse index  $m(P)$ equals the number of the sign changes
in the sequence
$$\mathcal{H}(P_3), \mathcal{H}(P_4),
\mathcal{H}(P_5),...,\mathcal{H}(P_n), $$ where the values
$\mathcal{H}(P_i)$ are already known by Theorem \ref{signHess}. \qed
\end{thm}

\begin{figure}\label{subconf}
\centering
\includegraphics[width=10 cm]{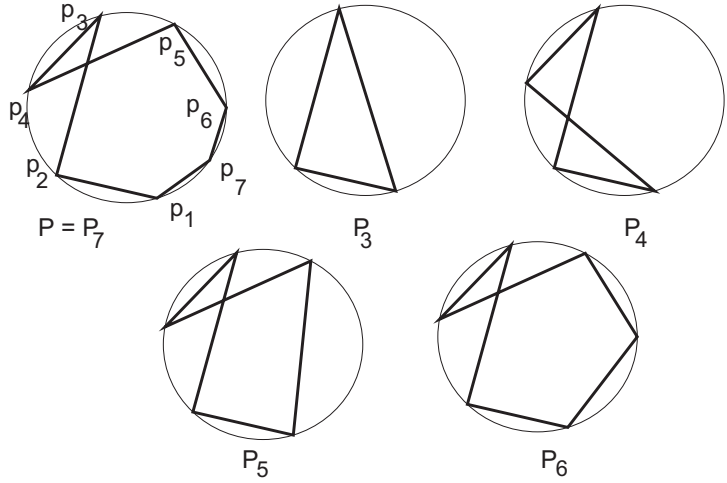}
\caption{Subconfigurations of $P$}
\end{figure}


\begin{thebibliography}{99}

\bibitem{cerf} J.Cerf, La stratification naturelle des espaces de fonctions
differentiables reelles et le theoreme de la pseudo-isotopie, Inst.
Hautes Et. Sci. Publ. Math. No. 39 (1970), 5-173.


\bibitem{eljikh}
E.Elerdashvili, M.Jibladze, G.Khimshiashvili, Cyclic configurations
of pentagon linkages, Bull. Georgian Nat. Acad. Sci. 2 (2008), No.4,
13-16.

\bibitem{FaS} M.Farber, D.Sch\"{u}tz, Homology of planar polygon
spaces, Geom. Dedicata 125 (2007),75-92.


\bibitem{GiNe}
C.Gibson, P.Newstead, On the geometry of the planar 4-bar mechanism,
Acta Applic. Math. 7 (1986), 113-135.



\bibitem{khi1}
G.Khimshiashvili, On configuration spaces of planar pentagons, Zap.
Nauch. Sem. S.-Peterb. Otdel. Mat. Inst. RAN 292(2002), 120-129.

\bibitem{khi2}
G.Khimshiashvili, Signature formulae and configuration spaces, J.
Math. Sci.  (accepted).

\bibitem{khipan}
G.Khimshiashvili, G. Panina, Cyclic polygons are critical points of
area.  Zap. Nauchn. Semin. POMI  360 (2008),  238--245.



\bibitem{rob}
D.Robbins, Areas of polygons inscribed in a circle, Discrete Comput.
Geom. 12(1994), 223-236.


\end{thebibliography}
\end{document}